# Extensionalism without Logicism: Ambrose and Extensional Logic

**Author:** Juan J. Colomina-Alminana

**Affiliation:** Northeastern University, Oakland Campus (Mills College)

**Address:** 5000 MacArthur Blvd., 132 Vera Long, Oakland, CA, 94613 – USA

**Email:** j.colomina-alminana@northeastern.edu

**Declaration:** The author declares that has no conflict of interest.

No data employed

This article adheres to all ethical standards

**Abstract:** Drawing primarily on her early work (1931-1934), I argue that Alice Ambrose's philosophical project at the time centers on preserving the rigor of extensional logic while rejecting the metaphysical and epistemological endorsements of logicism because of its commitment to the notion of material infinity. Positioning Ambrose as a transitional figure between Russell's formalism and the constructivist turn represented by Brouwer's and Weyl's intuitionism, I demonstrate how Ambrose offers a clever practice-oriented statement of finitist extensionalism. Employing only extensional methods—considering classes, relations, and propositions by reference to their members and truth-values instead of mental processes—, Ambrose reformulates an existential claim about $\pi$ as an explicit infinite disjunction of concrete instances insisting, against intensional projects, that such claims gain meaning only through a finite stopping rule that produces a witness.

* * *

This paper examines Alice Ambrose's nuanced engagement with extensionalism and, more specifically, Bertrand Russell's foundational work in logic and mathematics. While Ambrose often cites Russell's extensional methods as an example—his treatment of classes, relations, and propositions by reference to their members and truth-values—, she offers a sustained critique of his logicist ambitions, and particularly his endorsement of the materiality of infinity.

Bertrand Russell's contributions to logic and the foundations of mathematics are among the most influential of the twentieth century. His collaboration with Alfred North Whitehead on the *Principia Mathematica* (3 volumes, 1910–13; second edition, 1925-1927) poses that all mathematics could be derived from a small set of logical axioms and rules of inference. This project—known as logicism—was grounded in a formal, extensional calculus that treated classes, relations, and propositions as symbolic entities defined by their extensions. For Russell, the triumph of logicism lays on the realization that "all mathematics is symbolic logic" (Russell 1903: 44).

Alice Ambrose (1906-2001), American philosopher and formidable logician in her own right, engages deeply with extensionalism and, more concretely, Russell's philosophy of mathematics and logic in her early scholarly work. Specifically, in her 1932 unpublished dissertation "In Defense of an Extensional Logic" (defended at the University of Wisconsin-Madison), and in subsequent articles, Ambrose offers both praise and critique. She endorses Russell's extensional methods as the most rigorous and general framework for mathematical reasoning, particularly highlighting its advantage point versus intensional logical systems. At the same time, she rejects Russell's logicist claim that mathematics could be reduced to pure logic without remainder. She also challenges Russell's metaphysical commitments, particularly the materiality of infinity and his reliance on non-tautologous axioms like reducibility.

This paper, therefore, explores Ambrose's early extensionalist stance and, specifically, how it relates to Russell's logicist foundational program. I argue that her early defense of extensional logic is inseparable from her critique of logicism, and that her philosophical project preserves the methodological clarity of extensionalism while rejecting its metaphysical excesses. In other words, Ambrose defends an extensionalism without logicism, which is already present in her early writings but only becomes explicit when fully embracing a moderate finitism upon her arrival at Cambridge in the summer of 1932 (Cf. Ambrose 1935a and 1935b).[1]

Drawing primarily upon her early published and unpublished work (1931-1934), therefore, I demonstrate how Ambrose is a transitional figure during the decisive period of the foundational debates in logic and mathematics—one who bridges the gap between Russell's extensionalism and the emerging constructivist and finitist movements, particularly endorsed in the intuitionism demonstrated by Brouwer and Weyl. As proof of her groundbreaking middle ground, I show how Ambrose offers a clever practice-oriented statement of her finitist extensionalist method, one that feeds from the spirit of a proto-algorithmic solution. Relying only on extensional methods—considering classes, relations, and propositions by reference to their members and truth-values alone instead of their intensional meanings or mental processes—, Ambrose demonstrates that existence is epistemically tied to finite verification. Specifically, arguing against Black's intuitionistic solution to the so-called "π-7s problem," Ambrose reformulates an existential claim about π as an explicit infinite disjunction of concrete instances and insisted that such claims gain meaning only through a finite stopping rule that produces a witness, thus reducing existence claims to verifiable finite procedures and treating infinite disjunctions extensionally.

[1] There are a few recent articles addressing her work and her struggling relations with others during her time at Cambridge (1932-1935). Cf. Connell (2022), Chapman (2024), and Loner (2024).

**1. Ambrose on Russell's Extensionalism.**

Ambrose's admiration for Russell, and specifically his methodological approach, is clear from her early work, particularly that prior to her arrival in Cambridge as a postdoctoral student in the summer of 1932. As early as in 1929, she writes: "Mr. Russell plea[s] for rational scepticism. The burden of [Russell's] *Sceptical Essays* is that it is undesiderable to believe a proposition where there is no ground for supposing it true" (Ambrose 1929: 224). By endorsing extensional logic, I argue, that same skeptical spirit is also patently evident throughout her dissertation and early publications. Ambrose later writes: "A system of logic such as we find in the *Principia Mathematica* of A. N. Whitehead and Bertrand Russell pays its respects even more assiduously to premises than to consequences" (Ambrose 1932: 3). The emphasis on premises reflects the character of an extensional logic, where propositions are defined by their truth-values and classes by their members. Ambrose, therefore, aligns herself early on with this methodological approach, arguing that formal, extensional definitions yield clearer and more objective primitives for mathematical and logical inference and deduction.

Ambrose further endorses Russell's view that logic is prior to mathematics: "That such are the relative positions of logic and mathematics, I follow Bertrand Russell in maintaining against Brouwer and Weyl" (Ambrose 1932: 2). Ambrose sides with Russell against constructivists who ground mathematics on intuition and mental construction. For Ambrose, therefore, the formal, extensional calculus of the *Principia* provides a more rigorous and general foundation than that offered by any intuitionistic construct.

In her 1933 article "A Controversy in the Logic of Mathematics," Ambrose reiterates such an endorsement: "All of Russell's work in logic and the philosophy of mathematics is permeated with the assurance that logical foundations for mathematics can be laid down and that no arbitrary line between logic and mathematics can be drawn" (Ambrose 1933: 595). Ambrose also praises Russell's extensional infinite classes, which he considers "given all at once by the defining property of their members" (Russell 1916: 156), even though she understands that intuitionism has some advantage due to its endorsement of certain restrictions. She says: "The issue in fact rests on whether one can meaningfully define sets whose elements cannot be individually exhibited. That is to say, the validity of Zermelo's axiom of inclusion, which defines a sub-set by means of a property uniting the elements possessing it, is here at stake" (Ambrose 1933: 601). Ambrose argues for a clean and non-constructive definition of infinite sets, consistent with the extensional logic she also defends but that, paradoxically, as shall be demonstrated later, Russell's logicist program cannot resolve.

In her 1934 review of Max Black's *The Nature of Mathematics*, Ambrose goes further, directly criticizing Black's intuitionism, again highlighting the strengths of extensionalism: "All that needs to be known of these objects is stated in the axioms" (Ambrose 1934: 363). This reflects her commitment to a logical program that operates purely on the basis of extensions out of tautological principles—primarily truth-values and class membership—without invoking intensional or psychological notions. In fact, she had already used C. I. Lewis's 1929 fully intensional project in *Mind and the World-Order* as a punching bag (Cf. Ambrose 1931). Nevertheless, and this is important for what it follows, even though Ambrose endorses Russell's "sceptical rationalism" and, specifically, Russell's extensionalist program in logic, she will criticize his logicism and, in particular, the material interpretation of the infinity axiom.

**2. Russell's Logicism and Material Infinite.**

In defending extensional logic, Ambrose never loses sight of Russell's logicist ambition—to derive "all of mathematics" from formal logic alone. She persistently warns that this ambition depends on unproved, extra-logical (primarily ontological) assumptions, specifically those concerning infinity. The heart of her critique lies within two points. First, Russell's logicism cannot dispense with substantive axioms (notably the infinity and reducibility axioms) if classical logic and mathematics—in particular the theory of the infinite—is accommodated. Second, Russell's "all-at-once" infinite classes become a non-tautological postulate, contradicting his claim that pure logic contains "no new knowledge" beyond its primitive propositions (Cf. Russell & Whitehead 1910: v).

In fact, Russell's celebrated maxim—"all mathematics is symbolic logic"—rests on every mathematical proposition being reduced to a logical tautology or to a definition in logic's primitive vocabulary (Russell & Whithead 1910: xi). Yet, the *Principia*'s two editions (1910-1913 and 1925-27) quietly introduce three axioms that go beyond mere tautology (Cf. Linsky 2011). First, the infinity axiom asserts at least one class containing the empty set and closed under the "successor" operation. Without it, one cannot prove the existence of the infinite series of natural numbers. Second, the multiplicative axiom guarantees that any set with an infinite family of subclasses requires a choice-set containing exactly one element from each subclass. This underlies the construction of products of infinitely many sets. Third, the reducibility axiom collapses all higher-order predicates into extensionally equivalent first-order predicates, thereby restoring classes of classes as though they were classes of individuals. As Ambrose puts it, the *Principia* allows hidden axioms, what extensionalism never intended, making Russell's logicism untenable (see especially Ambrose 1932: Chapter 5, where she spells the technicality of the proofs).

Ambrose insists that these axioms are not logical tautologies but statements of fact, for which the *Principia* offers no formal proof of consistency. She observes: “Up to the introduction of non-tautologous principles, *Principia* fulfills the first requirement [of deriving all of mathematics from logic]. Proof that a proposition of the form ‘A is non-A’ can never be deduced, *Principia* does not attempt. And such a proof of consistency Hilbert is bound to provide” (Ambrose 1932: 74. See also Hilbert 1930: 245). In short, Russell’s logicist program achieves a formal unification with mathematics only by sneaking in substantive, unreduced assumptions—particularly infinity—contrary to logic adds “no new knowledge” beyond its formal structure.

In addition, Russell endorses an extensional view of the infinite, famously declaring that “classes which are infinite are given all at once by the defining property of their members, so that there is no question of ‘completion’ or of ‘successive synthesis!’” (Russell 1916: 156). Ambrose—even though sympathetic as she is to this extensional solution—warns that infinite as “given all at once” assumes the materiality of the infinity axiom, because one cannot otherwise guarantee that some class is not merely potentially infinite but actually (meaning, materially or physically) infinite. She writes: “Is a definition by inclusion equivalent to a definition by extension? The intuitionist answers no. His restriction is that definitions must contain the means of constructing the object defined. It is impossible to give a constructive definition, for example, of the set whose elements are the points of space, a set whose cardinal number is $2^{\aleph_0}$. For constructibility means essentially one-by-one exhibition. Existence in an infinite set is not a well-defined attribute, for an infinite process of logical addition may be requisite to call the indefinite element into being” (Ambrose 1933: 600). According to Ambrose, the extensional doctrine of immediate wholeness for infinite classes (Russell 1903: Chapter 43), dissolves into a hidden material axiom unrecoverable from logic’s primitive laws of identity, contradiction, and excluded middle. In other words, Ambrose argues that a genuinely infinite class cannot be justified purely through logical tautologies, highlighting a key tension between formal logic and mathematical intuition.

The reducibility axiom was Russell's ingenious—but ontologically hefty—patch to his theory of logical types. It avoids a vicious circle. Ambrose highlights its logicist dissonance though: "Russell brings in the axiom of reducibility which makes a characteristic of higher order equivalent in extension to one of lower order. That is, it is a scheme for reducing the order of a propositional function. Now the assumption of the existence of classes would make the axiom of reducibility unnecessary. But acceptance of the latter is a smaller assumption than that of the existence of classes; and most mathematicians will feel more comfortable, other things being equal, if 'entities are not multiplied beyond necessity.' The only logical reason for assuming classes is their use in reducing the order of a propositional function" (Ambrose 1932: 54). By forcing all higher-order predicates to be extensionally equivalent to some first-order predicate, Ambrose argues, the reducibility axiom makes every propositional function 'extensional'—yet does so only by asserting a new first-order predicate, a statement of fact wholly outside logic's tautologies. Far from closing the circle then, Ambrose understands that Russell's reducibility axiom leaves logicism unproven and adrift. In Ambrose's words, mathematics may be reduced to logic only at the cost of re-introducing mathematics as facts within logic.

Paradoxically, Russell's extensional clarity can be retained—classes by membership, propositional functions by extension, infinite sets "given all at once"—without embracing full logicism. By abstaining from asserting infinite classes (this is, by assuming what Ambrose calls the "no-classes" theory) and by eschewing the reducibility axiom, one arrives at a calculus of pure extensional relations. As she writes: "Russell felt convinced in 1906 that the plausibility on the surface of the simpler theories was truly only on the surface, and that though the no-classes theory was complicated, it was, nevertheless, preferable. Ordinary mathematics and most of the theory of the transfinite are even so constructible. Further, the theory is marked by one less assumption. Nothing is said as to the existence of classes, and such silence provides pleasing conformity with the ideal of Occam's razor" (Ambrose 1932: 50. Cf. Ambrose 1933). For Ambrose, the logic which treats classes, relations, and propositions purely by their extensions is free of all self-referential paradoxes provided we do not assume that every predicate determines a class or that every class is an individual. What remains is an extensional calculus of proof minus the material interpretation of infinity.

In sum, Ambrose's critique reveals that logicism overreaches. To secure Cantor's infinities and Russell's set theory, one must introduce extra-logical axioms that violate logic's parsimony. Ambrose, though, preserves extensional methods—truth-functions, classes by extension, "all at once" infinities—while dropping the logicist demand for a single unbroken derivation of mathematics from logic alone. Her solution is, therefore, to move toward finitism, defending an extensionalism without logicism.

**3. Type Theory, Reducibility, and the No-Classes Hypothesis.**

As Ambrose demonstrates, Russell's concern with the materiality of infinity emerges as a central tension within his program of deriving mathematics from purely logical principles. Russell's extensional logic relies on the assumption that infinite totalities can be treated as complete entities. This assumption, however, is paradoxical, most famously contradicting "the class of all classes that are not members of themselves." To resolve it, Russell introduces the theory of types in Appendix B of the *Principles*, aiming to prevent self-referential definitions by stratifying objects into hierarchical levels. As he writes: "What the complete solution of the difficulty may be, I have not succeeded in discovering; but as it affects the very foundations of reasoning, it is important to realize its existence" (Russell 1903: 525). Later, in the *Principia*, Russell and Whitehead formalize the so-called ramified theory of types, supplemented by the reducibility axiom, which controversially reintroduces impredicative definitions preserving classical analysis. This move intended to reconcile the logicist program with the material reality of infinite classes, though Russell himself remained uneasy about its philosophical legitimacy.

In his 1908 paper "Mathematical Logic as Based on the Theory of Types," Russell acknowledges: "In all the contradictions there is a common characteristic, which we may describe as self-reference or reflexiveness" (Russell 1908: 222). Ultimately, Russell's no-class theory—which treats classes as logical fictions—further avoids reifying infinite sets, preserving extensional logic while negating the ontological status of infinite collections. His work, thus, reflects a deep ambivalence: committed to logicism, yet haunted by the problematic materiality of infinity.

Russell's logicist program, therefore, was imperiled by paradoxes that arose when propositions or predicates range unrestrictedly over all objects—even over themselves. Ambrose identifies this impasse as the moment when Russell chose between two unpalatable alternatives. Russell was forced to either assume that every propositional function determines a class (thereby inviting his own "class of all classes" paradox) or abstain from committing to the existence of classes altogether. The means by which Russell navigates this impasse—his theory of types bolstered by the reducibility axiom— becomes the target of Ambrose's sustained critique and prompted her to explore the "no-classes" approach as a purer expression of extensional logic.

Besides his well-known correspondence with Frege, Russell's paradox emerges most famously when W. E. F. Burali-Forti drew attention to the "ordinal of all ordinals," a self-referential construction that demonstrates that any naïve set-theoretic system attempting to form the set of all ordinals leads to an immediate contradiction (Cf. Russell 1903: §301; Russell 1907; Jourdain 1907: 282).[2] In response, Russell proposes the logical types theory, in which each propositional function and class belong to a strictly ordered hierarchy. He writes: "Whatever involves all of a collection must not be one of the collection... A function cannot be its own argument, because the argument must be of a different logical type" (Russell & Whitehead 1910: 37-38). In this scheme, the natural number series of individuals, the classes of individuals, the classes of classes, and so forth, occupy distinct levels—preventing paradoxical self-reference. As he summarizes it: "We avoid the contradiction by arranging terms in a hierarchy of types, so that a term of one type can only be predicated of terms of a lower type" (Russell 1919: 212).

Nevertheless, Russell himself recognizes that the type-hierarchy creates unwelcome obstacles for mathematical constructions. In particular, the proof of key theorems—i.e., the least upper bound for sets of real numbers—requires treating a higher-type predicate as if it were of the same type as its instances. The proposed remedy is the reducibility axiom: “We assume that, corresponding to every function $\varphi$ of a given order, there is a predicative function $\psi$ which is extensionally equivalent to $\varphi$, i.e., such that $\varphi(x) \equiv \psi(x)$ for all values of x” (Russell & Whitehead 1910: 54).

[2] The contradiction discovered by Burali-Forti arises from the consideration of the class of all ordinals. This class must itself have an ordinal number, since it is well-ordered; but this ordinal number must be greater than any ordinal and therefore cannot be an ordinal. The classic version asks us to think of ordinals as a way to label positions: 0, 1, 2, 3, 4,… are all finite ordinals. Take now ω as the first infinite ordinal, and then come ω+1, ω+2, ω+3, ω+4,… and so on. Each ordinal, then, could be defined as the set of all smaller ordinals. Suppose now that we try to define a set Ω that contains all ordinal numbers. Since ordinals are well-ordered, Ω itself would be well-ordered as well, which means that Ω would have an order type—i.e., it would be an ordinal number. Nevertheless, this leads us to a paradoxical result, because if Ω is an ordinal and it contains all ordinals, then, it must contain itself. But, as established, no ordinal can be a member of itself. Therefore, if $\Omega \in \Omega$ leads to a contradiction, we must conclude that Ω cannot be a set. The only way to avoid this in a formal theory is to forbid the comprehension principle that lets us form "the set of all X satisfying…" without restriction in the first place. Russell's theory of logical types and modern Zermelo–Fraenkel set theory (with the axiom of separation), as Ambrose recognizes, both block the paradox by preventing one from treating "the class of all ordinals" as a legitimate set, by carefully restricting what kinds of sets can exist.

Ambrose trenchantly observes that this axiom is not a tautology but a factual assertion. She says: “[Russell] meets this difficulty by introducing the axiom of reducibility which makes any higher-order characteristic extensionally equivalent to some first-order characteristic. Thus the property of belonging to at least one of the classes of the aggregate is made equivalent to a property of the same order as those of classes… It seems to be a statement of fact and not a tautology” (Ambrose 1932: 83). According to Ambrose, therefore, because it introduces material existence claims—“there is a first-order $\psi(x)$ such that $\varphi(x) \equiv \psi(x)$”—the axiom violates Russell’s own dictum that logic adds no new knowledge except definitions.

Faced with the dual burdens of avoiding paradox and preserving extensionality, Russell eventually adopted what Ambrose calls the “no-classes theory”: “When we refuse to assert that there are classes, we must not be supposed to be asserting dogmatically that there are none. We are merely agnostic as regards them — like Laplace, we can say, ‘*je n’ai pas besoin de cette hypothèse*’” (Russell 1919: 184, as quoted in Ambrose 1932: 50). By abstaining from class-existence claims, rescuing Russell’s “sceptic rationality” methodology, Ambrose argues, no longer requires reducibility or infinity as assumptions—rather, one works syntactically, treating class-notation as an incomplete symbol whose only role is to abbreviate a formal definition.

Ambrose lauds this maneuver as the purest form of extensional logic. Avoiding any existential assumption about classes restores a fully extensional calculus, free of both the vicious circles of predicative self-reference and the *ad hoc* reducibility axiom. Negating the existence of the class as one means that the class is not to be treated as a full object in the sense one of its members is. Symbolically these are of lower type than the function. Thus, "the class of classes not members of themselves"-contradiction vanishes. And viewing classes as aggregates of members secures what the reducibility axiom would otherwise do—making every higher-order predicate extensionally equivalent to some first-order one—without having to postulate it. As Ambrose says: "The definitive feature of an extensional logic is its treatment of class-symbols as if they represented a plurality. Insofar as reference is made to a class's defining function, we refer to the (complementary) intensional aspect. But this aspect is not stressed. Classes are the values satisfying a propositional function. Ramsey makes more explicit the place of individuals in saying that a class is defined by a function which is a tautology with any member of the class as argument and a contradiction for anything else" (Ambrose 1932: 95).

Yet, as Ambrose recognizes, the no-classes approach demands sacrificing the infinity axiom. Without assuming any class existence, one cannot prove theorems requiring an actual infinite class—i.e., "there are infinitely many primes," "there are three consecutive 7s in π," etc. (Cf. Ambrose 1936b). In addition, Ambrose notes, one risks syntactic overreach. A purely syntactic calculus detaches entirely from the semantic notion of "there exists," potentially undermining the link to applied mathematics. In this light, Ambrose concludes that, while the "no-classes theory" represents the purest extensional stance, a balance is called for—one that preserves extensional definitions but admits limited existential axioms to support mainstream mathematics.

**4. Ambrose's Defense of Extensional Logic.**

After dismantling Russell's logicist pretensions and his material invocation of infinity, Ambrose turns to the rival account of inference championed by C. I. Lewis—specifically strict implication and the calculus of intensions—to show that extensional machinery remains both sufficient and superior for mathematical reasoning.

Lewis introduces a calculus of strict implication to supplement Boolean truth-values (true/false) with two further modalities: impossibility (~p) and necessary truth (□p). Material implication (p⊃q), in contrast, is simply defined by the extensional disjunction (p⊃q) = (~pvq), a statement about truth-values alone without any appeal to possibility or necessity. Lewis objects to material implication on the ground that: "A false proposition implies any proposition, true or false" and "A true proposition is implied by any proposition, true or false," are both wildly "useless" theorems, failing to capture the "proper" sense of inference (Lewis 1912: 523–524; Cf. Ambrose 1932: 36–37).[3]

Ambrose concedes that intensional nuance illumines aspects of our ordinary reasoning, but she insists that mathematics does not require Lewis's full panoply of strict modalities. She says: "[In material implication,] the only implication banned is 'a true proposition implies a false one.' Mr. Lewis suggests defining 'p implies q' as '~p^q,' where the conjunction sign symbolizes an intensional relation in which the rejection of one alternative binds one to accept the other—binds one in a way in which one is not bound in accepting just any proposition q on the basis of not-p. One can know the truth here while still ignorant of which alternative is true or whether both are true" (Ambrose 1932: 12).

[3] Ambrose writes: "Modern extensional theory of logic—of the Russellian sort—imposes somewhat of a strain on current meanings of inference and implication: Implication is such that a true proposition is implied by any proposition and a false proposition implies any proposition. This being the case, we could not uniquely draw out the consequences of a false premise, nor could we find what unique premises implied a given true statement" (Ambrose 1931: 369). This passage specifically criticizes the extensional (material) interpretation of implication, which Lewis rejects in favor of an intensional logical system. In that paper, Ambrose argues that material implication leads to counterintuitive results, such as any proposition implying a true one and a false proposition implying anything—highlighting why Lewis insists on intensional meaning as foundational for logical and epistemological clarity.

In Lewis's system the primitive notions are: 1. Propositions, 2. Negation, 3. Impossibility (~p), 4. Conjunction (∧), 5. Strict implication (p→q) and 6. Equivalence (=). From these, Lewis derives material implication and the standard Boolean operators, plus the four modalities: necessity (□p) and possibility (◊p), and their opposites (□~p) and (◊~p). Ambrose points out that three of these primitives can in fact be defined in purely extensional terms: (□p = ~◊~p), (◊p = ~□~p), and (p→q = ~◊p∧~q). The intensional operators, therefore, collapse into material ones when viewed extensionally. Ambrose writes: "The ideas of possibility, impossibility, and necessary truth are dispensed with as unnecessary complications to a working system of logic… They are seemingly not necessary to the utility of a logic—(and the uselessness of some of the theorems of Material Implication is the basis of [Lewis'] main protest.) I rather suspect that they reflect an underlying metaphysics. The world is such, he says [in (Lewis 1914)], that a difference exists between the false and the impossible, the true and the necessary, the real and the possible. And to such a world the system of material implication does not apply, for in its world these distinctions are erased. At the root of the distinction between the relations of strict and material implication lies the distinction between impossibility and simple falsity. In this fact, I suppose, rests its value for Lewis. For he insists that he has, by use of this category, given a 'proper' meaning to implication, without which proper meaning symbolic logic cannot be a criterion of valid inference" (Ambrose 1932: 112-113; Cf. Lewis 1918: 323-326).

Formal implication asserts a class of material implication for every value satisfaying the function. Ambrose shows that this extensional device, together with the rule of *modus ponens* for classes, reproduces every valid step of mathematical proof—no intensional modal machinery required. She emphasizes: "The crux of the issue is such terms as 'all' and 'there exists.' It is the consequence for logic of the intuitionist conception of 'all' and 'there exists' which has fanned the flames of controversy. The controversy centers about a second restriction, which the intuitionist imposes on the applicability of the law of excluded middle" (Ambrose 1933: 603). Inference in mathematics, Ambrose argues, never demands a modal premise beyond ~pvq and its universal closure, since formal implication plus extensional inference covers every case. Ambrose, therefore, corroborates what Norbert Wiener stated: "The certainty of mathematics rests on the 'uniformity' of our habits in handling ideas—on material implication and extensional definition, not on an a priori sense of possibility" (Weiner 1915: 571). If this is true, as Ambrose argues, by jettisoning intensional primitives altogether, one arrives at a purely extensional calculus based on the primitive notions of proposition, truth-value, negation, disjunction, and conjunction and where formal quantified implication P(x)R(x) is complemented with the material implication ~pvq. As Ambrose demonstrates, this minimal apparatus suffices for the deductive structure of all mathematics—set theory, number theory, analysis—once one accepts a limited postulate (i.e., the potentiality of infinity) but not the full logicist program of deriving the infinity axiom from pure logic nor its materiality (Cf. Ambrose 1932: Chapter 5). Ambrose, thus, concludes: "Pragmatic considerations… determine my approach to logic, and will appear in an important way, but their virtue will be extra-logical. […] It might seem—and of course this is contrary to the spirit of Brouwer—that the issue is one as to methodology, one which will be settled, not logically, but on pragmatic grounds" (Ambrose 1933: 600, 611).

**5. What Can Be Learned from Intuitionism?**

In the years immediately following her defense of an extensional logic, Alice Ambrose turns to the epistemological and methodological stakes of accepting actual infinities. Whereas in her prior work she exposes intensionalism, intuitionism, and logicism's extra-logical commitments, by her arrival in Cambridge in the summer of 1932 Ambrose had moved further toward a finitist outlook, fully embracing Brouwer's. The turn was then not unannounced, since she offers infinity animadversion ideas as soon as in 1931 and can be traced back to her dissertation. Nevertheless, she now explicitly insists that any assertion of an infinite mathematical object must be accompanied by a constructive method for exhibiting at least one of its members, siding with intuitionists such as Brouwer and Weyl. This section traces that shift to her later work, molding her adapted extensional program.

The intuitionist movement in the foundations of mathematics, led by L. E. J. Brouwer and later echoed by Hermann Weyl, challenges the classical conception of mathematics as a system grounded in logic and capable of manipulating infinite totalities. Their critiques—especially of the materiality of infinity—continue to resonate in philosophical debates about the nature of mathematical truth and the legitimacy of formal methods. Brouwer's intuitionism begins with a radical reorientation. According to Brouwer, mathematics is not a body of truths derived from logical axioms, but a mental activity rooted in the "Ur-intuition" of time (Brouwer 1907). From this intuition of temporal succession, Brouwer constructs the natural numbers and all valid mathematical operations. Logic, in his view, is not foundational but derivative—an abstraction of the mind from mathematical practice, not its source. As Ambrose summarizes: "Brouwer has answered this question [Is Logic basic?] with a flat no. Scorning any grounding of mathematics in logic, he has given primacy to mathematical intuition and forced logic to play second fiddle" (Ambrose 1933: 594).

Rejecting logic as a foundation is not merely rhetorical. Brouwer's constructive approach insists that mathematical objects must be built step by step, through finite mental acts. The infinite then is not a complete entity but a process—an unfolding that never reaches closure. Ambrose captures this in her exposition of Brouwer's "two-oneness" intuition.[4] She writes: "This intuition is the sense of a temporal unity in diversity, of a unit span of time in which inexhaustible division is possible… The intuition of two-oneness gives rise to the numbers one and two from which all finite ordinal numbers may be constructed and whose indefinite repetition creates the smallest infinite ordinal ω" (Ambrose 1933: 597). Here, the infinite is not material or actual—it is a method of construction, always incomplete. Brouwer's finitism thus denies the legitimacy of reasoning about infinite sets or totalities unless they can be generated through finite procedures.

Hermann Weyl, though initially aligned with Hilbert's formalism, later adopted a position closer to Brouwer's. He rejected impredicative definitions and the use of actual infinities in analysis, seeking instead a predicative and constructive foundation for mathematics. While Ambrose mentions Weyl only briefly, her remark is telling: "As one might expect, Weyl [(1929)] rejects the axiom of reducibility" (Ambrose 1932: 83; Ambrose 1933: 607). This rejection aligns with Weyl's broader critique of classical logic and set theory. For Weyl, the continuum is not a complete set of real numbers but a "medium of free Becoming"—a process that cannot be captured by static, extensional definitions. Like Brouwer, Weyl insists that mathematics must remain within the bounds of what can be constructed by finite means.

Ambrose's engagement with intuitionism, therefore, is not merely descriptive. She uses both Brouwer's and Weyl's views to highlight deep problems in the classical conception of mathematics and logic. She challenges that infinite sets or logical systems are materially real—i.e., complete objects. The intuitionist critique exposes classical mathematics' reliance on symbolic manipulation divorced from intuitive content. By foregrounding Brouwer's constructive method then, Ambrose shows that many classical proofs—especially those involving the law of excluded middle or non-constructive existence claims—lack genuine mathematical meaning. They assert the existence of entities that cannot be constructed and, thus, violate the intuitionist standard of legitimacy. Moreover, Ambrose's treatment of Brouwer's "two-oneness" intuition underscores the philosophical depth of intuitionism. Mathematics, in this view, is not a static edifice but a dynamic unfolding of mental acts. The infinite remains, then, not a thing but a direction—a horizon that guides construction without ever being reached.

[4] This is the name Ambrose gives to the notion of "*Zweieinigkeit*," which we now know as "two-ity" after the translation from the original written in Dutch by A. Heyting and E. Beth (Cf. Brouwer 1975). Brouwer (1907: 5) reads: "Mathematics is the languageless construction originating in the basic intuition of the falling apart of a life moment into two qualitatively different things, one giving way to the other, yet held together in the mind." And a few years later, he again describes mathematics as "the free creation of the mind, founded on the perception of the move of time" (Brouwer 1949: 483). For an early presentation of Brouwer's views, see Brouwer (1913).

What can, then, be learned from intuitionism? As Ambrose reminds us, first, one can learn that mathematics is not merely formal (like Russell's logicism assumed) but also experiential. The insistence on finite construction forces us to confront the epistemological basis of mathematical knowledge; this is, as Ambrose argues, what one can actually know, and how. Second, intuitionism challenges the dominance of logic in foundational debates. Brouwer's rejection of logic as primary invites a reconsideration of the role of intuition, temporality, and mental activity in mathematical reasoning. Ambrose takes this as fundamental when she accepts that, perhaps, at some point, pragmatic elements (extra-logical, as she often refers to them) have a role establishing logical principles. Third, criticizing the materiality of infinity is relevant to set theory, computability, and the philosophy of mathematics. As Ambrose reminds us, the intuitionist perspective warns against reifying the infinite—treating it as a completed object—when it is, in fact, a conceptual tool grounded in finite cognition.

In sum, Brouwer and Weyl, at least as interpreted by Ambrose, offer a compelling vision of mathematics as a human activity rooted in intuition and constrained by the limits of construction. Their rejection of actual infinity and impredicative reasoning provoke reflection on the nature and scope of mathematical truth. She has no other choice. She abandons all logicist pretensions and reformulates extensionalism's scope to better account for the foundations of logic and mathematics, incorporating a more practical stance into her methodological account.

**6. Ambrose's Balanced Extensional Logic.**

By her arrival in Cambridge during the summer of 1932, Ambrose had fully embraced a finitist turn, which materializes with the publication of Ambrose (1935a) and Ambrose (1935b), which famously demonstrates that her approach is at odds which the more extreme finitist account defended by Wittgenstein. Ambrose (1936a), for example, begins by concisely stating the intuitionist requirement that commanded her approval: "A necessary condition for the existence of a number is the presence of a method of determining it" (Ambrose 1936a: 117). This simple dictum—echoing Brouwer's emphasis on constructive exhibition—represents a decisive departure from earlier acceptance of Russell's "given all at once" infinities. Here, existence is no longer a vacuous instantiation of an extensional predicate but a requirement of finite demonstrability.

A specifically relevant instance of this new methodological extensionalism is also in Ambrose (1936a), where she revises the "π–7s problem." Black's *Claims of Intuitionism* had famously introduced the "π–7" rule: "Define $N$ as the real number in (0, 1) whose $k^{th}$ decimal digit is '1' if and only if $k$ consecutive 7s occur at the start of the decimal expansion of π, and '0' otherwise" (Black 1936: 91). Black argues that one must, in principle, scan all the digits of π up to the first 7-run of length $k$ in order to determine the $k^{th}$ digit of $N$. Ambrose corrects Black's logical conclusion, which according to her rests on an amendable error: "He mistakenly holds that it is necessary to examine *the whole* series 1, 4, 1, 5, 9,… in order to write *any* term of $N$. This is only necessary if 0 is written, and the point is that in advance there is no way of determining whether this is necessary or whether we can write 1. Black could, however, correct his example by requiring that *k* consecutive 7s occur exactly *n* times, for a fixed *n*, or by requiring that *k* consecutive 7s occur exactly *k* times" (Ambrose 1936a: 117). With this revision then, Ambrose demonstrates that finitism can be fully compatible with a broad class of "non-constructive" definitions, departing from the Wittgensteinian stand and pushing for an amended, balanced extensional logic without the issues of neither logicism nor intuitionism. One only needs a uniform recipe guaranteeing termination in finite time. Once that is satisfied, the existential claim attains a clear extensional meaning in accordance with logic's primitive rules.

Although Russell had championed the view that infinite classes are "given all at once" (Russell 1916: 156), which Ambrose once endorsed, by 1935-1936 she insists that even Cantorian infinities may be approached one-by-one. Remember that she writes: "It is impossible to give a constructive definition, for example, of the set whose elements are the points of space, a set whose cardinal number is $2^{\aleph_0}$. For constructibility means essentially one-by-one exhibition. Existence in an infinite set is not a well-defined attribute, for an infinite process of logical addition may be requisite to call the indefinite element into being" (Ambrose 1933: 600). Reframing existence in terms of finitary stopping rules,[5] Ambrose forges a middle path between Russell's material infinite and Brouwer's outright negation of classical inferences. Her finitism preserves the extensional logic's algebra of classes and truth-functions while securing an epistemic foothold: To assert any existential claim such as ∃(x)φ(x) is to provide a finite "witness" (a) such that φ(a) holds.

[5] Ambrose seems to be here anticipating the concept of calculability defended by Church and further solidified in the notion of a Turing machine. After all, she is highly aware of Church's work on number theory, as demonstrated by her work notes of the time (See Alice Ambrose Lazerowitz: Papers, GBR/0012/MS Add.9938/3/1, hosted at the University of Cambridge Libraries Special Collections).

Deeping even further, Ambrose (1936b) is a concrete application of her balanced extensional methodology. Ambrose begins by showing that infinite disjunction can work as an extension in certain cases. She argues that the question "Are there three consecutive 7s in π?" presupposes that π's expansion is a complete totality. Yet, since π is non-repeating and infinite, one can never fully materialize its digits. Therefore, how can such a question be meaningful or answerable without invoking some metaphysical or material assumption? This is the core at the "π-7 problem": How to resolve questions about infinity without appealing to materiality?

Ambrose proposes a balanced extensional solution in proto-algorithmic spirit. She argues that meaning arises from potential mathematical proof, not from actual material verification. Even if we cannot observe all digits of π, one understands what it would mean to verify the presence of three consecutive 7s in a finite allows meaning retention, even in the absence of material completion. In other words, Ambrose's proto-algorithmic solution is a philosophical maneuver. She shifts the burden of meaning from material instantiation to logical possibility, which aligns with her broader finitist leanings. Specifically, for Ambrose, the assertion 'There are three 7s in the $n^{\text{th}}$, $n+1^{\text{st}}$, and $n+2^{\text{nd}}$ places in the decimal development of π' is in effect one of an infinite disjunction—one disjunct for each *n*—namely:

[place Figure 1 here]

In which:

- $\pi_n$ refers to the $n^{\text{th}}$ digit in a sequence—in this case, the decimal expansion of π.
- ∧ refers to the Logical AND, the conjunction.

- $\sum$ over $n$=1 refers to the logical OR, the disjunction, and is equivalent to the application of the infinite disjunction over all digits $n$ of the sequence $n \geq 1$.

The whole expression, therefore, asserts: "There exists some position $n$ such that three consecutive terms—$\pi_n$, $\pi_{n+1}$, and $\pi_{n+2}$—are all equal to 7." According to Ambrose then, employing the disjunction tool as restricted method allows proving that there is at least one place in the sequence where three 7s appear in a row in the extension of π (Ambrose 1936b: 507-509).

This is exactly the extensional move in a logical set-theoretic context, or "in number theory," as she calls it (Ambrose 1936b: 511). This is a compact way to express a search for a triple repetition of a value. One can simply treat "Are there three consecutive 7s?" as the extension of the predicate "three consecutive 7s occur starting at position *n*," ranging over all *n*. This way, no intensional meaning is required. As she stresses: To know "what three consecutive 7s look like" is nothing else that to know "how to write down three consecutive 7s by some finite calculation" (Ambrose 1936b: 507)—never to invoke an intensional or meaningful connection between π and these digits. This, as mentioned above, keeps intact the extensional logic's dictum that "classes which are infinite are given all at once by the defining property of their members" (Russell 1916: 156), without carrying over the contradiction of infinity as a complete physical object.

Ambrose's endorsement of finitism against claims for the materiality or actuality of infinity shows that one cannot meaningfully pose "Are there three 7s at the $n^{th}$, $n + 1^{st}$, $n + 2^{nd}$ places, or not?... [because] we do not in one sense know what a method for proving one of an endless disjunction is like; and the answer which is the formal negative of this… we do not know in any sense" (Ambrose 1936b: 509). One only knows a finite stopping rule that halts whenever a run of three 7s appears. Nevertheless, there is no analogous finite test that confirms their non-occurrence across the entire expansion. In other words, we can simply assert this claim and take the 'Or–Not' part as a *non sequitur* without appealing to any constructive test because, as she mentions, "The phrase 'proof that p' in mathematics… has no significant opposite" (Ambrose 1936b: 511).

One interesting ontological consequence of Ambrose balanced position is that she rescues the extensional program without accepting the burden of the logicist claims. For Ambrose, existence requires finite witness. An existential claim (such as "$\exists n$…") must come with a finite rule that exhibits a particular witness $n$, a concrete point of stoppage within the sequence—there is no need to appeal to a Platonic or material infinite $n$. In her dissertation she had already faulted Russell's infinity axiom as a "non-tautologous principle" he "never proves" but simply "asserts" (Ambrose 1932: 71–72). Requiring finite stoppage recasts the existence of runs of 7s in π into purely finitist—and thus epistemically and algorithmically secure—terms. As she insists: "To write down three 7s is by no means to write down the result of the calculation of π… A mistake made in answering a mathematical question by doing a calculation… consists in writing down figures not in conformity with those rules… Since one has not set up rules for a different calculation, one is in a sense doing no calculation" (Ambrose 1936b: 512).

On the one hand, Ambrose's π-digit example implements her thesis that extensional logic handles infinite phenomena by treating them as "given all at once" sets of truth-value assignments. On the other, she applies that extensionalism to a specific infinite disjunction about π's digits, augmenting it with a finitist requirement: That any ∃-claim must be backed by a concrete, finite stopping rule. Ambrose, therefore, implicitly rejects the materiality of infinity as a primitive while preserving the clarity of the extensional method.

## 7. Conclusion.

This study traced Alice Ambrose's evolving engagement with extensionalism, demonstrating how she consistently champions its formal, extensional machinery—classes defined purely by their members, propositional functions handled solely by truth-values, and infinite sets "given all at once"—yet resolutely rejects the logicist ambitions that smuggle substantive axioms into logic's tautologous framework.

From her 1931 published critical notice of Lewis's book onward, Ambrose affirmed that an extensional, formal apparatus offers the most rigorous, ontology-minimal basis for mathematics. By treating all relations as extensions—classes of ordered pairs or sets of truth-value assignments—extensionalism achieves a unity of form that Ambrose found methodologically unbeatable.

Ambrose methodically criticizes logicism for its hidden metaphysical and epistemic assumptions, keeping Russell on check for that. The logicist program's insistence that "all mathematics is symbolic logic" (Russell 1903: vii. Cf. Russell & Whitehead 1910: 3) falters on the necessity of the infinity and reducibility axioms. Ambrose shows that these axioms convert logicism into a metamathematical enterprise, thereby abandoning its original ideal of a self-sufficient logic.

As this paper has demonstrated, Ambrose endorses a reform of material infinity. By the mid-1930s, Ambrose explicitly embraces a principled finitism. She says: "Existence in an infinite set is not a well-defined attribute, for an infinite process of logical addition may be requisite to call the indefinite element into being" (Ambrose 1933: 603. Cf. Ambrose 1935a). Her correction of Black's response to the "π–7 problem" (Ambrose 1936a) further exemplifies this stance. Existence, therefore, is an extensional predicate fulfilled by some exhibited member—no actual infinite need to be assumed as a primitive beyond pragmatic warrant.

Ambrose, therefore, defends a balanced extensional program for the foundations of mathematics and logic, one that preserves the power of extensional calculi—truth-functions, quantification over classes, and the "given-all-at-once" approach to infinite collections—while relaying on mathematical intuition. This allows Ambrose to offer a proto-algorithmic proof to the "π-7 problem" based on finite witnesses and stopping points. Ambrose thus abandons the logicist drive to extract every mathematical theorem from logic's primitive tautologies and self-imposes an epistemic check—any ∃-claim must be backed by a finite method of construction, thereby aligning mathematical practice with human cognitive and mechanical capabilities while distancing from the idea that mathematics is a purely formal endeavor.

Ambrose's balanced extensional program, hence, anticipates important currents in post-1936 logic and computer science. For example, the "no-classes" hypothesis dovetails with modern NF-style set theories that reject a global Russell paradox by restraining comprehension without elaborate type hierarchies, as endorsed by, for example, modern Type-Free Set Theory. Her finitist requirement that ∃-claims be accompanied by witnesses, mirrors witness-extraction in proof-theoretic semantics, where classical theorems are given constructive content. Revision of the infinite decimal-rule (the π-7s extension problem above) presages computable real numbers, each presented by a Turing-decidable stopping rule—an embodiment of Ambrose's finitist extensionalism in algorithmic form.

Once again, Ambrose's middle grounds methodology between intuitionism and logicism attempts of founding logic and mathematics documents why an extensional program should replace material infinity with potential possibilities. Ambrose's finitist axioms yield new hierarchies of predicative versus impredicative power, which obliges to re-examining the precise proof-theoretic strength of any logical system. This, though, we leave for another day.

**Appendix. Interpreting the Ambrose Formula as a Formal Rule**

One can treat the formula in Figure 1 as a logical condition that defines a rule for selecting positions in the decimal expansion of π. Let's formalize it and express it as a computable function. This expression says: For some index ($n$), the digit at position ($n$) in π is 7, and so are the digits at positions ($n$ + 1) and (n + 2) (since ($n$ + 2 = n+1)). Therefore, such a condition simplifies to: $\pi_n = 7 \wedge \pi_{n+1} = 7$. This defines a trigger rule: We are looking for consecutive 7s in the decimal expansion of π. This formalizes in the following formal rule:

Let $(d_i)_{i \in N}$ be the digits of π after the decimal point.

Define the set of trigger positions:

$T = \{i \in N \mid d_i = 7 \wedge d_{i+1} = 7$

For each $i \in T$, define a block $B(i)$ of the next 7 digits after the second 7:

$B(i) = (d_{i+2}, d_{i+3}, \ldots, d_{i+8})$

Therefore, the Ambrose-style 7s-extension sequence (what we can dub the Ambrose Formula) is:

E = Concatenate all $B(i)$ for $i \in T$

This rule exemplifies Ambrose's view that mathematical meaning is not intuitive but rule-governed, nor intensional but extensional. Since the expression defines a publicly followable condition, the extension is determined by use, and not by mental insight. Therefore, this proto-algorithm is a model of meaning through mechanical procedure, not intuition—exactly the kind of example Ambrose uses to challenge Black's position, and the one that ended antagonizing Wittgenstein's views on the foundations of mathematics.

## Acknowledgements

Deleted for peer-review purposes.

Figure 1.

$$\sum_{n=1} \pi_n = 7 \wedge \pi_{n+1} = 7 \wedge \pi_{n+2} = 7]$$